\documentclass[12pt,leqno]{article}
\usepackage{amsthm, amsfonts}
\begin{document}
\newtheorem{theo}{Theorem}[section]
\newtheorem{cor}{Corollary}[section]
\newtheorem{lem}{Lemma}[section]
\newtheorem{pro}{Proposition}[section]
\newtheorem{app}{}[section]
\def\RR{{\cal R}}     
\def\pl{{{\Bbb P}^1}} 
\def\C{{\Bbb C}}      
\def\Pf{\C[t]}         
\def\proj{{\Bbb P}^2} 
\def\O{{\cal O}}      
\def\s{{\cal S}}
\def\M{{\cal M}}
\def\N{{\Bbb N}}
\def\Z{{\Bbb Z}}
\def\Q{{\Bbb Q}}
\def\R{{\Bbb R}}
\def\H{{\cal H}}
\def\U{{\cal U}}
\def\A{{\cal A}}
\newcommand{\RE}[2]{\Omega^{#1}_{#2}}
\newcommand{\REg}[2]{\overline{\Omega}^{#1}_{#2}(*D)}
\newcommand{\deri}[2]{\nabla{#2}{#1}}
\newcommand{\tih}[1]{\H^{#1}}
\newcommand{\REa}[2]{\tilde{\Omega}^{#1}_{#2}}

\begin{center}
{\LARGE\bf Relative Cohomology with Respect to a Lefschetz Pencil
\\}
\vspace{.25in} {\large {\sc Hossein Movasati}\footnote{Supported
by IMPA-Brazil, IPM-Iran
\\
Math. classification: 14F05, 14F43 \\
Technische Universit\"at Darmstadt 
Fachbereich Mathematik,
Schlo\ss gartenstr. 7, 
64289 Darmstadt, Germany   
Email: {\tt movasati@mathematik.tu-darmstadt.de}, Fax: 0049-6151-162747 
}}
\end{center}
\begin{abstract}
Let $M$ be a complex projective manifold of dimension
$n+1$ and $f$ a meromorphic function on $M$ obtained by a generic
pencil of hyperplane sections of $M$. The $n$-th  cohomology vector bundle of 
$f_0=f|_{M-\RR}$, where $\RR$ is the set of 
indeterminacy points of $f$, is defined on the set of
 regular values of $f_0$ and we 
have the usual Gauss-Manin connection on it. Following Brieskorn's methods in
\cite{bri}, we extend the $n$-th cohomology 
vector bundle of $f_0$ and the associated Gauss-Manin connection to $\pl$ by
means of differential forms.
The new connection  turns out to be meromorphic on the  
critical values of $f_0$. We prove that the meromorphic 
global sections of the vector bundle with poles of arbitrary order at
$\infty\in\pl$ is isomorphic to the Brieskorn module of $f$ in a natural way,
and so the Brieskorn module in this case is a free $\Pf$-module of rank
$\beta_n$, where $\Pf$ is the ring of polynomials in $t$ and
$\beta_n$ is the dimension of $n$-th cohomology group of a regular fiber of
$f_0$. 
\end{abstract}
\setcounter{section}{-1}
\section{Introduction}
The algebraic description of the monodromy of a germ of 
an isolated  singularity $f:(\C^{n+1},0)\rightarrow(\C,0)$ was done by 
E. Brieskorn in \cite{bri}. In this article he considers the Milnor fibration 
associated to $f$ and then the $n$-th cohomology vector bundle $\H$ of $f$ over a punctured neighborhood of $0\in\C$ and
the associated Gauss-Manin connection. Then he constructs three extension of
this vector bundle (the sheaf of its holomorphic sections) $\H,{'\H},{''\H}$  
by means of holomorphic forms in $(\C^{n+1},0)$. He constructs them up to
torsions  which may appear in the
stalk over zero
but later M. Sebastiani in \cite{seb} proves that there is no torsion and
so Brieskorn's extension is complete. By a slight modification of his argument we can obtain a meromorphic connection 
$\nabla:V\rightarrow \Omega^1_{\C,0}(k0)\otimes_{\O_{\C,0}}V$
which is the Gauss-Manin connection of the Milnor fibration in $(\C,0)-\{0\}$.
Here $V$ stands for one of $\H,{'\H}$ and ${''\H}$, 
$k$ is the smallest number with this property that the multiplication by $f^k$ induces the zero
map in the Jacobi algebra $\frac{\O_{\C^{n+1},0}}{<f_{x_i}\mid i=1,\ldots,n+1>}$ of $f$ and
$\Omega^1_{\C,0}(k0)$ is the sheaf of meromorphic 1-forms in $(\C,0)$ 
with poles of order at most $k$ at $0$. 
The stalk  of $\H,{'\H},{''\H}$ over $0\in(\C,0)$, namely $H, {'H}, {''H}$,
are called Brieskorn modules and they are very useful objects in 
singularity theory. They are freely generated $\O_{\C,0}$-modules of rank $\mu$,
where $\mu$ is the Milnor number of $f$.

We consider a projective manifold $M$ of dimension $n+1$
 and we intersect it by a pencil of hyperplanes (see for instance Lamotke's
article \cite{lam} 
for the definitions concerning a pencil).
We can define a holomorphic map $f_0$ in $M-\RR$, where $\RR$ is the 
intersection of the axis $A$ of the pencil with $M$, such that $\overline{f_0^{-1}(t)}$'s are hyperplane sections. 
We assume that $A$ intersects $M$ transversally and $f_0$ has only isolated
singularities. Now $f_0$ is a $C^{\infty}$ fiber bundle over $\pl-C$, where $C$
is the set of critical values of $f_0$. Therefore we have the cohomology fiber 
bundle over $\pl-C$ and the associated Gauss-Manin connection. 

In the above context we will generalize Brieskorn methods as follows: We
make a blow up $\pi:V\rightarrow M$ along $\RR$ and we obtain our extensions by means of meromorphic forms in $V$ with poles of arbitrary order along $\pi^{-1}(\RR)$.
In this way we obtain three analytic sheaves ${\H}^n,{'\H}^n,{''\H}^n$ on $\pl$
and connections 
$\nabla:W\rightarrow\Omega^1_{\pl}(\tilde C)\otimes_{\O_{\pl}}W$.
Here $W$ stands for one of ${\H}^n,{'\H}^n$ and ${''\H}^n$, 
$\tilde C$ is a divisor in $\pl$ with support at the critical values of
$f_0$ and $\Omega^1_{\pl}(\tilde C)$ is the sheaf of meromorphic 1-forms in 
$\pl$ with pole divisor less than or equal $\tilde C$.

Our main result in this article
is that $W$ is a locally free sheaf of rank $\beta_n$
 over $\pl$, where $\beta_n$ is the dimension of the $n$-th cohomology of a 
regular fiber of $f_0$, and there is a canonical
isomorphism between $W\mid_{\pl-C}$ and the $n$-th cohomology vector bundle of $f_0$ over 
$\pl-C$
and $\nabla$ is the Gauss-Manin connection by this isomorphism.
 Then we introduce global Brieskorn modules $H^n,{'H}^n, {''H}^n$ in our context.
They are $\Pf$-modules, where $\Pf$ is the ring of polynomials in $t$, 
and we prove that there is a 
$\Pf$-module isomorphism between $H^n$ (resp. ${'H}^n$ and ${''H}^n$) and the module of global meromorphic sections
of ${\H}^n$ (resp. ${'\H}^n$ and ${''\H}^n$) with poles of arbitrary order
at $\infty\in\pl$.

In the lower dimensions $i<n$ there is no vanishing cycle and the monodromy around a critical value
is identity. It is not difficult to see that in this dimension the cohomology
of the critical fiber gives us the desired extension and the Gauss-Manin connection on the $i$-th cohomology vector bundle is holomorphic even in
the critical value. This implies that the obtained vector bundle is trivial. However, we construct this extension by means of meromorphic
forms.

The notion of global Brieskorn modules has been recently considered by many 
people, see for instance C. Sabbah, A. Dimca, M. Saito and P. Bonnet's 
works  \cite{sab1}, \cite{disa}, \cite{bon}. As an immediate consequence of our last result
we prove that $H^n$ and ${'H}^n$ and ${''H}^n$ are freely generated $\Pf$-modules of rank $\beta_n$. This result was already known by 
C. Sabbah in \cite{sab1}. 
In the context of differential equations (n=1) ${'H}^1$ appears in 
the works of G.S.
 Petrov for polynomials of the type $y^2+P(x)$ in $\C^2$ and is called Petrov
module by L. Gavrilov.
(see \cite{gav2} p. 572). Recently some applications of this module
in differential equations
have been introduced by the author of these lines in \cite{mov3}.

{\bf Brieskorn module ${'H}$ in differential equations:} Consider the case
$M=\proj$ and $f=\frac{F}{G}$, where $F$ and $G$ are two
polynomials of the same degree in an affine coordinate $\C^2$ of
$\proj$. Assume that $F=0$ intersects $G=0$ transversally and the
critical points of  $\frac{F}{G}$ are non-degenerate with distinct
images. Consider the foliation
\begin{equation}
\label{sextasanta}
{\cal F}_\epsilon: df+\epsilon.\omega=0
\end{equation}
where $\omega$ is a meromorphic 1-form in $\proj$ with poles of arbitrary order
along $G=0$.
Let $\{\delta_t\}_{t\in(\C,0)}$ be a continuous family of vanishing cycles.
We call $h(t):=\int_{\delta_t}\omega$ an Abelian integral.
If $h(t)\not\equiv 0$ then the cycle $\delta_{t_0}$ persists in
being cycle after this deformation if and only if $h(t_0)= 0$.
Therefore the study of the number of limit cycles appearing from
$\delta_t$'s after the deformation ${\cal F}_t$  leads to the
study of the zeros of Abelian integrals. In \cite{mov1} 
(see also \cite{mov2}) it is shown that if $h(t)\equiv0$ then $\omega=Pdf+dQ$,
 where $P$ and $Q$ are two meromorphic functions in $\proj$ with
poles of arbitrary order along $G=0$. Therefore $\omega=0$ in ${'H}^1$ and ${'H}^1$
represents the space of deformations ~(\ref{sextasanta}) for
which the birth of limit cycles can be studied by Abelian
integrals. 

Now let us explain the structure of this article.
In \S ~\ref{27may02} we have explained in details the extension of 
the cohomology vector bundles of $f_0$ to the critical values of $f_0$ 
and the associated  Gauss-Manin connection to a meromorphic connection.
Theorem ~\ref{khodaya} which is the central result in this article is stated 
there together with Theorem  ~\ref{bazasheg}.  The reader who is interested only on the construction and the main results is invited to read only this section.
Theorem  ~\ref{bazasheg} is proved in \S ~\ref{mehdi}. 
The terminology and propositions in Appendix ~\ref{17mar01} are used in this
section.  
 \S ~\ref{mehdi2} is devoted to the proof of Theorem ~\ref{khodaya}. 
The proof can be considered as a kind of variational Atiyah-Hodge theorem, therefore it is recommended 
to the reader to know the proof of this theorem stated in \cite{nar}.

Perhaps three appendices for this article is too many, but in each of them
we have obtained some partial results which we need them in this article and
I did not find them in the literature.
In Appendix ~\ref{17mar01} we have listed some necessary concepts and theorems in complex geometry. The first result is ~\ref{A9} 
which is a kind of Kodaira vanishing theorem for direct limit of coherent 
 sheaves. ~\ref{8abr01} is the main result in this appendix.
It is a kind of variational Kodaira vanishing theorem  and is frequently used 
in this article. After doing a blow-up in the indeterminacy locus of our pencil
we obtain a holomorphic map $g:V\rightarrow \pl$ and a divisor $A$ in $V$ in such a way that the intersection of $A$ with a regular fiber is a positive divisor in that fiber.
 ~\ref{8abr01} claims that $R^ig_*\s(*A)=0, i>0$, where $\s$ is a coherent sheaf in $V$ and $\s(*A)$ is the sheaf of meromorphic section of $\s$ with poles
of arbitrary order along $A$. Note that $g$ has critical fibers.
 In Appendix ~\ref{alman} we list some information about the
topology of the fibers of $f_0$. Any kind of singularities can appear in our
pencil therefore we had to combine some technics of \cite{lam} and \cite{arn}
to obtain ~\ref{20ago01}, the main result of this appendix. In particular we prove that a distinguished basis of vanishing cycles in the singularities with a same value must be linearly independent.   
Let us consider the restriction map from a global Brieskorn module of a pencil 
to a local Brieskorn module of a singularity of the pencil. It is believed that this map is surjective but I was not able to prove this fact.
In Appendix ~\ref{gudiny} we prove ~\ref{9ago01} which
says that the local Brieskorn module divided by the image of the mentioned 
map is a vector space of finite dimension. This easily implies that the image of the mentioned map is a freely generated $\O_{\C,0}$-module, the statement of 
~\ref{kalisch}. 
    
I would like to thank my teachers   S. Shahshahani, C. Camacho, A.
Lins Neto, and P. Sad for their
support and interest. I would like to thank also C. Hertling, S. Lu  
and Y. Holla for
many useful conversations. The first draft of this article is obtained in
IMPA-Brazil and IPM-Iran.  The final version is obtained in MPIM-Bonn. 
Here I acknowledge my sincere thanks to all these institutes 
for hospitality. 
\section{Extension of the cohomology vector bundle and the Gauss-Manin
connection}
\label{27may02}
Let $M$ be a complex projective manifold of dimension
$n+1$, $\{M_t\}_{t\in\pl}$ a pencil of hyperplane sections of $M$
and $f$ the meromorphic function on $M$ whose level sets are
$M_t$'s. Suppose that the axis of the pencil
intersects $M$ transversally (see \cite{lam}). This implies that the set of
indeterminacy points $\RR$ of $f$ is a smooth submanifold of $M$
of codimension two and every two $M_t$'s intersect each other transversally 
in $\RR$.
Define $f_0:=f\mid_{M-\RR}$ and $L_t:=M_t-\RR=f_0^{-1}(t)$.
We assume that the critical points of $f_0$ are isolated and 
we denote by $C=\{c_1, c_2,c_3,\ldots,c_r\}$ the set of critical
values of $f_0$. Note
that a critical fiber $M_{c_j}$ may have more than one critical
point. Put
\[
\beta_i=dim(H^i(L_t,\C)),\ t\in\pl-C,\ 0\leq i\leq n
\]
$f_0$ is a $C^\infty$ fibration over $\pl-C$ (see for instance
\cite{lam}) and so $\beta_i$ is independent of $t$.

The set $\tilde{H}^i=\cup_{t\in \pl-C}H^i(L_t,\C)$ has a natural structure 
of a complex manifold and the natural projection $\tilde{H}^i\rightarrow \pl-C$
is a holomorphic vector bundle map which is called the $i$-th cohomology 
vector bundle. 
Let $\C_{M-\RR}$ be the sheaf of constant functions in $M-\RR$
and  $R^if_{0*}{\C_{M-\RR}}$ be the $i$-th direct image of the
sheaf $\C_{M-\RR}$ (see \cite{gra2} and Appendix ~\ref{17mar01}). 
Any element of
$R^if_{0*}{\C_{M-\RR}}(U)$, $U$ being an open set in $\pl-C$, is a
holomorphic section of the cohomology vector bundle map and is
called a constant section. It is easy to verify that
\[
\O(\tilde{H}^i)\cong R^if_{0*}{\C_{M-\RR}}\otimes_{\C}\O_{\pl-C}
\]
where $\O(\tilde{H}^i)$ denotes the sheaf of holomorphic
sections of $\tilde{H}^i$. We
define the sheaf
\[
\tilde{\H}^i= R^if_{0*}{\C_{M-\RR}}\otimes_{\C}\O_{\pl-C},\ 0\leq\ i\leq n
\]
The Gauss-Manin connection on $\tilde{\H}^i$ is given by
\begin{equation}
\label{22may}
\deri{}{}:\tilde{\H}^i \rightarrow \Omega^1_{\pl-C}  \otimes_{\O_{\pl-C}}\tilde{\H}^i 
\end{equation}
\[
\deri{}{(g\otimes c)}=dg\otimes c,\ c\in R^if_{0*}{\C_{M-\RR}}(U),\ g\in \O_{\pl-C}(U)
\]
where $U$ is an open set in $\pl-C$. Now we have the problem of extension of $\tilde{\H}^i$ to a locally free sheaf on $\pl$ and $\nabla$ to a (meromorphic)
connection (with possible poles in $C$) defined in the extended sheaf. We could define
$\tilde{\H}^i$ and $\nabla$ in the whole $\pl$. 
In  ~\ref{20ago01} Appendix ~\ref{alman} we have proved that 
$\tilde{\H}^i, \ 0\leq i<n$ is a locally free sheaf of rank $\beta_i$ in $\pl$
so in this case there is no serious problem. 
But in the case $i=n$ we have the notion of vanishing cycle in a critical
point of $f_0$ and so the definition ~(\ref{22may}) does not
give us the desired extension. Using Brieskorn's ideas we are going to 
construct
some extensions  and reconstruct $\tilde \H^i$'s by means of meromorphic forms.
 
Let $\pi:V \rightarrow M$ be the blow-up along $\RR$ (see \cite{lam}). $V$ is a smooth manifold and
$g=\pi\circ f$ is a well-defined holomorphic function in $V$. We have
\begin{equation}
\label{22}
A:=\pi^{-1}(\RR)\cong\pl\times \RR,\ V_t:=g^{-1}(t)=\pi^{-1}(M_t),\ 
\ L_t\cong V_t-A
\end{equation}
Each fiber $V_t$ intersects $A$ transversally in $\{t\}\times\RR$.
For $U\subset \pl$ we define 
\begin{equation}
\label{11}
V_U:=g^{-1}(U),\ L_U:=V_U-A
\end{equation}
Let $\tilde{\Omega}^i$ be the sheaf of holomorphic $i$-forms in $V$ and
$\tilde{\Omega}^i(*A)$ be the sheaf of its meromorphic sections with poles of arbitrary order along $A$. Let also $\Omega^i$ be  the direct image by $g$ of
$\tilde{\Omega}^i(*A)$, i.e. $\Omega^i=g_*\tilde{\Omega}^i(*A)$. 
The following sheaf is well-defined
\[
\Omega^{i}_{V/\pl}=\frac{\Omega^i}{\Omega^1_{\pl}\wedge{\Omega}^{i-1}}
\]
We have the following long sequence:
\begin{equation}
\label{26may02}
\ 0\stackrel{d^{-1}}{\rightarrow}
\Omega^{0}_{V/\pl}\stackrel{d^0}{\rightarrow}
\Omega^{1}_{V/\pl}\stackrel{d^1}{\rightarrow}\cdots\stackrel{d^{n-1}}
{\rightarrow}
\Omega^{n}_{V/\pl}\stackrel{d^{n}}{\rightarrow} \Omega^{n+1}_{V/\pl}
\end{equation}
We define
\[
\tih{i}=\frac{ker d^i}{Im d^{i-1}}
\]
Note that $\Omega^{i}_{V/\pl}$ are $\O_{\pl}$-module and 
the differential operators $d^i$'s are $\O_\pl$ linear and so the cohomology sheaves $\tih{i}$'s are  $\O_\pl$-modules. We also define the sheaves:
\[
{'\H}^n=\frac{\Omega^n}{\Omega^1_{\pl}\wedge \Omega^{n-1}+d\Omega^{n-1}}
,\ 
{''\H}^n=\frac{\Omega^{n+1}}{\Omega^1_{\pl}\wedge d\Omega^n}
\]
From now on the index $i$ stands for $0,1,2,\ldots,n,'n,''n$. For instance 
if $i='n$ then $\tih{i}={'\H}^n$ and $\beta_i=\beta_n$.
To construct connections on $\H^i$'s we need the following lemmas.
The key of the proof in both lemmas is ~\ref{8abr01} in Appendix ~\ref{17mar01}.
Let $c\in\pl$, $U$ a small open disk around $c$
and $t$ a regular holomorphic function in $U$.
\begin{lem}
(generalized de Rham lemma)
\label{dohafte}
An element $\omega\in\Omega^i(U), i\leq n$ is of the form $dt\wedge \eta,\ \eta\in\Omega^{i-1}(U)$ if and only if $dt\wedge\omega=0$.
\end{lem}
\begin{proof}
It is enough to prove that if
 $dt\wedge \omega=0$ then $\omega$ is of the form $dt\wedge\eta,\ \eta\in\Omega^{i-1}(U)$.
By de Rham lemma (see \cite{bri}, p. 110) we can write
\[
\omega=dt\wedge \eta_\alpha,
\eta_\alpha\in\tilde{\Omega}^{i-1}(*A)(U_\alpha)
\]
where $\{U_\alpha\}_{\alpha\in I}$ is an open covering of $V_U$.
Now
$\{\eta_\alpha-\eta_\beta\}_{\alpha,\beta\in I}$ is an element of
$H^1(V_U, \s(*A))$, where
\[
\s=Ker(\tilde{\Omega}^{i-1}\stackrel{dt\wedge
}{\rightarrow}\tilde{\Omega}^{i})\
\]
By ~\ref{8abr01}  $H^1(V_U, \s(*A))=0$
and so we can find $\eta'_\alpha\in \s(*A)(U_\alpha)$ such that
$\eta_\alpha-\eta_\beta=\eta'_\alpha-\eta'_\beta$. The
$(i-1)$-form $\eta\mid_{U_\alpha}=\eta_\alpha-\eta'_\alpha$
satisfies $\omega=dt\wedge \eta$ and is the desired $(i-1)$-form.
\end{proof}
\begin{lem}
\label{aykhoda}
$\Omega^{n+1}_{V/\pl}$ is a discrete sheaf with support at $C$. The stalk of
$\Omega^{n+1}_{V/\pl}$ over $c_j\in C$ is a vector space of dimension
$\mu_{c_j}$, where $\mu_{c_j}$ is the sum of Milnor numbers of the 
singularities of $f$ within $L_{c_j}$. 
In particular there is a natural number $k_j$ such that
$(f-c_j)^{k_j}\Omega^{n+1}_{V/\pl}$ is zero in $c_j$. 
\end{lem}
\begin{proof}
Put $\mu_c=0$ if $c$ is not a critical value of $f$.
We know that
\[
(\tilde\Omega^{n+1}_{V/\pl})_x:=
(\frac{\tilde\Omega^{n+1}}{g^*\Omega^1_{\pl}\wedge \tilde\Omega^n })_x=
\]
\begin{equation}
\label{24may02} 
              \left\{ \begin{array}{ll}
           \C^{\mu_x} & x \hbox{ is a critical point of $g$ 
          with the Milnor number } \mu_x \\
          0  & \hbox{ otherwise} \\
               \end{array}
    \right.
\end{equation}
By  ~\ref{8abr01} $H^1(V_U,g^{*}(\Omega_\pl^1)\wedge\tilde{\Omega}^{n}(*A))=0$ and
so
\[
\Omega^{n+1}_{V/\pl}(U)= \frac {H^0(V_U,
\tilde{\Omega}^{n+1}(*A))}{H^0(V_U, 
g^{*}(\Omega_\pl^1)\wedge\tilde{\Omega}^{n}(*A))}= 
H^0(V_U,\frac{\tilde{\Omega}^{n+1}(*A)}
{g^{*}(\Omega_\pl^1)\wedge\tilde{\Omega}^{n}(*A)})=
\]
\[
H^0(V_U,\frac{\tilde{\Omega}^{n+1}}
{g^{*}(\Omega_\pl^1)\wedge\tilde{\Omega}^{n}}(*A))
\]
Since  ~(\ref{24may02}) is a discrete sheaf with support at $V-A$, we conclude
that 
\[
\Omega^{n+1}_{V/\pl}(*A)(U)=\oplus_{x\in V_c}(\tilde\Omega^{n+1}_{V/\pl})_x=\C^{\mu_{c}}
\]
where $x$ runs through all critical points of $f$ within $V_c$.
\end{proof}
We put $k_j$ the minimum number with the property in Lemma ~\ref{aykhoda}.
For every
critical point $p$ in the fiber $V_{c_j}$ there exists a natural
number $k_p$ depending only on the type of the critical point
$p$ such that $(g-c)^{k_p}(\REa{n+1}{V/\pl})$ is zero in $p$ 
(see \cite{bri}, p. 110 and 125). Choose always the minimum $k_p$. 
We have $k_j=\max_p\{k_p\}$, where $p$ runs through all critical points of $f$ within $L_{c_j}$. 

Consider the sheaf $\tih{i},1\leq i\leq n-1$. Let
$[\omega]\in\tih{i}(U)$. We can write $d\omega=dt\wedge\eta,
\eta\in \Omega^i(U)$. We have $dt\wedge d\eta=0$ and so by Lemma
~\ref{dohafte} we have $d\eta=dt\wedge \eta'$ for some $\eta'\in
\Omega^{i}(U)$. Therefore we can define the following connection:
\[
\deri{}{}:\tih{i}\rightarrow\Omega^1_{\pl} \otimes_{\O_{\pl}}\tih{i} 
\]
\[
\deri{}{[\omega]}= dt\otimes[\eta],\  d\omega=dt\wedge \eta
\]
Now we are going to define the connections
 \[
\deri{}{}:\tih{i}\rightarrow\Omega^1_{\pl}(\tilde C)\otimes_{\O_{\pl}}\tih{i},\
i=n,'n,''n 
\]
where $\tilde C=\sum k_jc_j$ and $\Omega^1_{\pl}(\tilde C)$ is the sheaf of meromorphic
sections of $\tih{n}$ with a pole of order at most $k_j$ at $c_j$.
The reader is referred to \cite{esn} for more informations about meromorphic
connections on sheaves on $\pl$.

Let us define $p(t)$ in $U$ as follows: $p(t)=1$ if $U\subset \pl-C$ and $p(t)=(t-c_j)^{k_j}$ if $U$ is an open disk around
$c_j$. 
Let $[\omega]\in \H^n(U)$. We can write $d\omega=dt\wedge\eta,
\eta\in \Omega^n(U)$. Since $d\eta$ may not be of the form 
$dt \wedge \eta',\ \eta'\in \Omega^n(U)$, we had to multiply 
$d\eta$ by $p(t)$ and therefore by Lemma ~\ref{aykhoda} we have:
$$
d(p(t)\eta)=p(t)d\eta+p'(t)dt\wedge\eta=
dt \wedge \eta'+ p'(t)dt\wedge\eta,\ \eta'\in \Omega^n(U)
$$
Therefore we can define
\[
\deri{}{[\omega]}=\frac{dt}{p(t)}\otimes[\eta']
,\  \eta'=p(t)\eta,\ d\omega=dt\wedge\eta
\]
In a similar way  for $[\omega]\in {'\H}^n(U)$ 
\[
\deri{}{[\omega]}=\frac{dt}{p(t)}\otimes[\eta], \ p(t)d\omega=dt\wedge\eta
\]
and for $[\omega]\in {''\H}^n(U)$ we have 
\[
\deri{}{[\omega]}=\frac{dt}{p(t)}\otimes [d\eta], \ p(t)\omega=dt\wedge\eta
\]
It is not difficult to see that these definitions are well-defined and 
do not depend on the choice of the coordinate $t$ and the choice of $\omega$ in the
class $[\omega]$. From now on we use 
the notation $\omega$ instead of $[\omega]$. The main theorem of this
article is:
\begin{theo}
\label{yeksal}
\label{khodaya}
$\tih{i}, i=0,1,\ldots,n,'n, ''n$ is a locally free sheaf of rank $\beta_i$ on
$\pl$. The natural map $\tih{i}\rightarrow {\tilde\H}^i$ in $\pl-C$ which is obtained by the restriction of differential forms to the fibers of $g$  
induces an isomorphism between $(\tih{i},\nabla)$ and $({\tilde\H}^i,\nabla)$.
\end{theo}
We have used the convention $\tilde\H^i=\tilde\H^n$ for $i='n,''n$.
If we consider only one fiber $V_t$ then by Atiyah-Hodge theorem
(see ~\ref{hodge}) we know
that meromorphic differential forms in $V_t$ with poles of arbitrary order
along $A\cap V_t$ give us the cohomology groups of $V_t-A$. This shows 
that the above theorem in $\pl-C$ is a natural statement. Main difficulty
in the proof of the above theorem lies in the critical
values of $f$. To prove it we will have to look more precisely to 
the proof of Atiyah-Hodge theorem stated in \cite{nar}.

Now we can look at $\H^i$ as a vector bundle. In the case $i<n$ the obtained
connection is holomorphic in $\pl$. This implies that the vector bundle 
$\H^i, i<n$ is a trivial bundle. We have already expected this fact.  

Note that the above extensions are not necessarily logarithmic. All logarithmic extension 
to $C$ are described in \cite{anbo}, p. 89 and also in \cite{esvi}, \cite{her}
 for arbitrary dimension of the base space. 
  The extensions introduced above have a peculiar property which we are 
going to explain below:

Choose  $p=\infty\in\pl-C$. 
This implies that $D:=M_\infty$ is smooth. Let 
$t$ be an affine coordinate of $\C=\pl-\{p\}$ and  $\overline{\Omega}^{i}(*D)$ be the set of meromorphic $i$-forms in $M$ with 
poles of arbitrary order along
 $D$. 
 It is a $\Pf$-module in a trivial way
\[
p(t).\omega=p(f)\omega,\ \omega\in \overline{\Omega}^{i}(*D),\  p(t)\in\Pf
\]
where $\Pf$ is the ring of polynomials in $t$. 
The set of  relative meromorphic $i$-forms with poles of arbitrary order
 along $D$ is defined
as follows:
\[
\REg{i}{V/\pl}=\frac{\overline{\Omega}^i(*D)}{df\wedge\overline{\Omega}^{i-1}(*D)}
\]
The differential operator
\[
d^i: \REg{i}{V/\pl}\rightarrow\REg{i+1}{V/\pl}
\]
\[
\omega_1\rightarrow d\omega_1
\]
is well-defined and $\Pf$-linear. Now we have the complex of
relative meromorphic  forms $(\REg{*}{V/\pl},d^*)$ and so we can
form the cohomology groups
\[
H^i=H^i(\REg{*}{V/\pl}, d^*)=\frac{Ker(d^i)}{Im(d^{i-1})},\
d^{-1}=0
\]
$H^i$ is called the $i$-th relative cohomology of $M$ with respect to $f$.
It is easy to see that $H^{0}=\Pf$. In dimension $n$ there are 
two other useful $\Pf$-modules
\[
{'H}^n=\frac{\overline{\Omega}^n(*D)}{df\wedge \overline{\Omega}^{n-1}(*D)+
d\overline{\Omega}^{n-1}(*D)},\ 
{''H}^n=\frac{\overline{\Omega}^{n+1}(*D)}{df\wedge d\overline{\Omega}^n(*D)}
\]
The $\Pf$-modules $H^n,{H'}^n, {H''}^n$ were
introduced by Brieskorn in \cite{bri} for a germ of a holomorphic
function $f:(\C^{n+1},0)\rightarrow (\C,0)$ and are also called
Global Brieskorn modules (see Appendix ~\ref{gudiny}).

We consider the $\Pf$-morphism
\begin{equation}
\label{lu1}
H^i\rightarrow H^0(\pl, \tih{i}(*p)),\ i=0,1,2,\ldots,n,'n,''n
\end{equation}
and the morphism of vector spaces
\begin{equation}
\label{lu2}
\REg{n+1}{V/\pl}\rightarrow H^0(\pl, \RE{n+1}{V/\pl}(*p))
\end{equation}
obtained by restrictions to the fibers of $g$, where $\tih{i}(*p)$ denotes 
the sheaf of meromorphic sections of $\tih{i}$ with
poles of arbitrary order at $p$.
The next main result in this article is:
 \begin{theo}
\label{bazasheg}
The morphisms ~(\ref{lu1}) and ~(\ref{lu2}) are isomorphisms. 
\end{theo}
\begin{cor}
Keeping the notations used above, we have
\label{18mar01}
\begin{enumerate}
\item
$\REg{n+1}{V/\pl}$ is a vector space of dimension $\mu$, where
$\mu$ is the sum of local Milnor numbers of $f$;
\item
  $H^{i},\ i=0,1,\ldots,n,'n,''n$ is a free  $\Pf$ module of rank $\beta_i$;
\end{enumerate}
\end{cor}
\begin{proof}
By Theorem ~\ref{bazasheg} the first statement is trivial. For the second one
it is enough to
prove that $H^0(\pl, \tih{i}(*p))$ is a free $\Pf$-module of rank $\beta_i$.
 The sheaf $\tih{i}$ is a locally free sheaf of rank $\beta_i$ over $\pl$.
  By Birkhoff-Grothendieck decomposition theorem (see \cite{gra3}) there exist integers
  $n_1,n_2,\ldots,n_{\beta_i}$ (uniquely determined up to a permutation) such that
\[
\tih{i}\cong \O(n_1 p)\oplus\O(n_2
p)\oplus\cdots\oplus\O(n_{\beta_i} p)
\]
Meromorphic sections of $\O(n_i p)$ with poles of arbitrary order at
$p$ is a $\Pf$-module of rank one and therefore $H^0(\pl,
\tih{i}(*p))$ is a free $\Pf$-module of rank $\beta_i$. 
\end{proof}

If $i<n$ then the vector bundle $\H^i$ is trivial and so the numbers
$n_1,n_2,\ldots,n_{\beta_i}$ are zero. 
It would be interesting, if one tries to understand the nature of the numbers
$n_1,n_2,\ldots,n_{\beta_i}$ by some numerical invariants of the manifold
and the singularities of $f_0$ in the case $i=n,'n,''n$.

 The above corollary generalizes Brieskorn and Sebastiani's results
 in \cite{bri} and \cite{seb} in the local case $f:(\C^{n+1},0)\rightarrow (\C,0)$. When I finished this article I was informed that similar results
to Corollary ~\ref{18mar01} 
are obtained by C. Sabbah \cite{sab1} in the context of 
algebraic geometry and therefore in view of Serre's GAGA principal it is true
in the context of analytic geometry.
\section{Proof of Theorem ~\ref{bazasheg}}
\label{mehdi}
\newcommand{\relz}[1]{\Omega_\pl^1\wedge\Omega^{#1}}
\newcommand{\REaa}[2]{\Omega^{#1}_{#2}}
Recall the notations  introduced in ~(\ref{11}) and ~(\ref{22}).
We will also use notations introduced in Appendix ~\ref{17mar01}.
Let $\tilde{\Omega}^i(kA)$ be the sheaf of meromorphic $i$-forms
in $V$ with poles of order $\leq k$ along $A$ and
$\Omega^i(k)=g_*\tilde{\Omega}^i(kA)$.
By Grauert direct image theorem  
$\Omega^i(k)$ is a coherent sheaf and so we have a direct system of coherent sheaves $\{\Omega^i(k)\}_k$ (the map $\Omega^i(k)\rightarrow \Omega^i(k+1)$ is the 
inclusion). The sheaf $\Omega^i=\lim_{k\to \infty} \Omega^i(k)$
is the direct image of $\tilde{\Omega}^i(*A)$.
 The quotient
direct system of sheaves
\[
\{\RE{i}{V/\pl}(k)\}_k=\frac{\{\Omega^i(k)\}_k}
{\{\Omega^1_{\pl}\wedge\Omega^{i-1}(k)\}_k}
\]
is called the direct system of sheaves of relative meromorphic
$i$-forms. We can easily check that the differential operator
\[
d^i: \{\RE{i}{V/\pl}(k)\}_k\rightarrow\{\RE{i+1}{V/\pl}(k+1)\}_k,\ 
d^i(\omega)=d\omega
\]
is well-defined and $\O_\pl$-linear. Now we have the following, not necessarily exact, sequence
\[
\{ \RE{0}{V/\pl}(k)\}_k
\stackrel{d^0}{\rightarrow}\{ \RE{1}{V/\pl}(k+1)\}_k\stackrel{d^1}{\rightarrow}\cdots
\stackrel{d^{i-1}}{\rightarrow}\{\RE{i}{V/\pl}(k+i)\}_k\stackrel{d^{i}}{\rightarrow}\cdots
\]
We have the complex of relative meromorphic  $i$-forms
 $(\{\RE{*}{V/\pl}(k+*)\}_k, d^*)$, and so we can form the direct system of cohomology sheaves
\[
\{\tih{i}(k)\}_k=H^i(\{\RE{*}{V/\pl}(k+*)\}_k, d^*)=
\frac{\{Ker(d^i)\}_k}
{\{Im(d^{i-1})\}_k},\  i\geq 0,\ d^{-1}=0
\]
The differential operators $d^i$'s are $\O_\pl$ linear and so $\{\tih{i}(k)\}_k$ is a direct system of $\O_\pl$-module sheaves. We also define the sheaves:
\[
\tih{i}=\lim_{k\to\infty} \tih{i}(k)
\]
\[
\{{'\H}^n(k)\}_k=
\frac{\{\Omega^n(k)\}_k}{\{\relz{n-1}(k)+d\Omega^{n-1}(k-1)\}_k},
\ {'\H}^n=\lim_{k\to\infty} {'\H}^n(k)
\]
\[
{''\H}^n(k)=\frac{\{\Omega^{n+1}(k)\}_k}{\{dt\wedge d\Omega^n(k-1)\}_k},\ {''\H}^n=\lim_{k\to\infty} {''\H}^n(k)
\]
It is not clear at all that the sheaf $\tih{i}$
 is coherent. By properties ~\ref{A1} and ~\ref{A2} in Appendix ~\ref{17mar01}
 it is not difficult to see that this way of definition of $\H^i$'s coincide with
the one introduced in the first section.
\begin{theo}
The sheaves ${\H}^i(k),\ i=0,1,\ldots,n,'n,''n, k=0,1,2,\ldots$ are coherent.
\end{theo}
\begin{proof}
According to Grauert direct image theorem ( Theorem ~\ref{gdit} in Appendix ~\ref{17mar01})
$\Omega^i(k)$ is coherent for all finite integer number $k$.  Now the 
coherence of
$\tih{i}(k), k\in \N$ can be obtained from the following facts: 
Let $\s$ and $\s'$ be two coherent sheaves on a complex manifold $P$. 
If $\s\subset\s'$ then $\frac{\s}{\s'}$ is coherent. If $i:\s\rightarrow\s'$ is a $\O_P$-linear map then kernel and image of $i$ are coherent sheaves (see \cite{gra2} p. 236-237).
\end{proof}
Now let us prove Theorem ~\ref{bazasheg}.  We only prove the isomorphism 
~(\ref{lu1}). The proof of the second is similar.
First we observe that
\[
\tih{i}(*p)=\frac{Ker(d^i)}{Im(d^{i-1})}(*p)
\stackrel{~\ref{A4}}{=}\frac{Ker(d^{i})(*p)}{Im(d^{i-1})(*p)}
\]
The number over an equality means the corresponding property in
Appendix ~\ref{17mar01}. $Im(d^i)$ is a direct limit of coherent
sheaves and so by ~\ref{A9} we have $H^1(\pl,Im(d^{i-1})(*p))=0$ and
\[
H^0(\pl, \tih{i}(*p))= H^0(\pl,
\frac{Ker(d^i)(*p)}{Im(d^{i-1})(*p)})\stackrel{~\ref{A9}}{=}
\frac{H^0(\pl, Ker(d^i)(*p))}{H^0(\pl, Im(d^{i-1})(*p))}
\]
Let
\[
\cdots\rightarrow H^0(\pl,
\RE{i-1}{V/\pl}(*p))\stackrel{d^{i-1}_p}{\rightarrow} H^0(\pl,
\RE{i}{V/\pl}(*p))\stackrel{d^i_p}{\rightarrow} H^0(\pl,
\RE{i+1}{V/\pl}(*p))\rightarrow\cdots
\]
be obtained from ~(\ref{26may02}). By ~\ref{A5} and ~\ref{shayad} we have
\[
Ker(d^i_p)=H^0(\pl, Ker(d^i)(*p)),\  Im(d^{i-1}_p)=H^0(\pl, Im(d^{i-1})(*p))
\]
$\relz{i-1}$ is a direct limit of coherent
sheaves, by ~\ref{A9}  $H^1(\pl, (\relz{i-1})(*p))=0$ and
\[
H^0(\pl, \RE{i}{V/\pl}(*p))=H^0(\pl,
\frac{\Omega^i(*p)}{(\relz{i-1})(*p)}) \stackrel{~\ref{A9}}{=}
\frac{H^0(\pl, \Omega^i(*p))}{H^0(\pl,
(\relz{i-1})(*p))}\stackrel{~\ref{A7}}{=}
\]
\[
\frac{H^0(\pl, \Omega^i(*p))}{dt\wedge H^0(\pl, \Omega^{i-1}(*p)))}
\]
where $t$ is the chart map for $\pl-\{p\}$ (it has a pole of order
one in $p$). We have
\[
\Omega^i(*p)=(g_*\tilde{\Omega}^i(*A))(*p)\stackrel{~\ref{A10.5}}{=}g_*\tilde{\Omega}^i(*A+
*V_p), \ \pi^{-1}(D)=A+ V_p
\]
By blow down along $A$ we can see easily that $H^0(\pl
,\Omega^i(*p))\cong \overline{\Omega}^i(*D)$.
\qed
\section{Proof of Theorem ~\ref{khodaya}}
\label{mehdi2}
\newcommand{\diss}[2]{R^{#1}g_*#2}
\newcommand{\REab}[2]{{\cal E}^{#1}_{#2}}
\newcommand{\ses}[3]{0\rightarrow #1 \rightarrow #2 \rightarrow #3 
\rightarrow 0}
\def\CC{{\cal C}}
The arguments of this section can be considered as a variational
Atiyah-Hodge theorem (see ~\ref{hodge} in Appendix ~\ref{17mar01}). 
It is highly recommended to the reader to
know the proof of Atiyah-Hodge theorem stated in \cite{nar}. First we will
prove the assertion of Theorem ~\ref{khodaya} for $i=0,1,\cdots,n-1,'n$. 
The same statements for $i=n,''n$ follows directly. 

We have constructed $\tih{i}$ by means of meromorphic forms in
$V$ with poles of arbitrary order along $A$. The following lemma enables us to
reconstruct it by means of holomorphic forms in $V-A$.

\begin{lem}
\label{diagram} 
Let $g:X\rightarrow Y$ be a continuous map
between paracompact Hausdorff spaces and suppose that two
complexes $\A$ and $\A'$ of Abelian sheaves over $X$ are given
together with mappings $h$ such that the diagram
\begin{equation}
\begin{array}{ccccccccc}
0 & \stackrel{{d'}^{-1}}{\rightarrow} & \A_0'& \stackrel{{d'}^0}{\rightarrow}  & \A_1' & \stackrel{{d'}^1}{\rightarrow}  &\A_2' & \stackrel{{d'}^2}{\rightarrow}  & \cdots \\
 &  & h_0\downarrow &  & h_1\downarrow & &h_2\downarrow & & \\
0 & \stackrel{{d}^{-1}}{\rightarrow}  & \A_0& \stackrel{{d}^0}{\rightarrow}   & \A_1 &  \stackrel{{d}^1}{\rightarrow}  &\A_2 & \stackrel{{d}^2}{\rightarrow}   & \cdots \\
\end{array}
\end{equation}
is commutative. (The rows are not supposed to be exact, but we
have $d\circ d=0$ and $d'\circ d'=0$). Suppose further that
\[
\diss{i}{\A_k}=0, \diss{i}{\A'_k}=0,\ \forall i\geq 1, k\geq 0
\]
and for $k\geq 0$ $h$ induces isomorphisms of cohomology sheaves
\begin{equation}
\label{18mai01} \frac{Ker ({d'}^k)}{Im({d'}^{k-1})}\rightarrow
\frac{Ker ({d}^k)}{Im({d}^{k-1})}
\end{equation}
Then $h$ induces isomorphisms
\begin{equation}
\frac{Ker ({d_*'}^k)}{Im({d_*'}^{k-1})}\rightarrow \frac{Ker
(d_*^k)}{Im(d_*^{k-1})}
\end{equation}
for all $k\geq 0$, where $d_*$ and $d_*'$ define the sequences
\begin{equation}
\begin{array}{ccccccccc}
0 & \stackrel{{d_*'}^{-1}}{\rightarrow} & g_*\A_0'&
\stackrel{{d_*'}^0}{\rightarrow}  & g_*\A_1' &
\stackrel{{d_*'}^1}{\rightarrow}  &g_*\A_2' &
\stackrel{{d_*'}^2}{\rightarrow}  &
\cdots \\
 &  & \downarrow &  & \downarrow & &\downarrow & & \\
0 & \stackrel{d_*^{-1}}{\rightarrow}  & g_*\A_0&
\stackrel{d_*^1}{\rightarrow}
 & g_*\A_1 &  \stackrel{d_*^2}{\rightarrow}  &g_*\A_2 & \stackrel{d_*^2}{\rightarrow}   &
  \cdots \\
\end{array}
\end{equation}
\end{lem}
\begin{proof}
We have just rewritten Theorem 6.5 of \cite{nar} in another
form.
\end{proof}
Let ${\cal E}^i$  be the sheaf on $V$, which is defined by the
presheaf that to every open subset $U$ of $V$ associated the
modules of holomorphic $i$-forms in $U-A$. Let also
\[
\REab{i}{V/\pl}=\frac{{\cal E}^i}{g^*\Omega^1_{\pl}\wedge{\cal
E}^{i-1}}
\]
Let $U$ be a small open disk in $\pl$. Since $L_U$ is a Stein
manifold (see ~\ref{20ago01}) and the restriction of any Stein
covering (see \cite{gra2}) of $V_U$ to $L_U$ is again a Stein
covering, by Cartan's B theorem we have
\[
\diss{j}{\REab{i}{V/\pl}}=0,\ j>0
\]
We have the following long sequence:
\begin{equation}
\label{tamam} \A :=\ 0\rightarrow
\REab{0}{V/\pl}\stackrel{d^0}{\rightarrow}
\REab{1}{V/\pl}\stackrel{d^1}{\rightarrow}\cdots\stackrel{d^{n-1}}{\rightarrow}
\REab{n}{V/\pl}\rightarrow 0
\end{equation}
 Recall that $\tilde{\Omega}^i(*A)$ is the sheaf of meromorphic $i$-forms in
$V$ with poles of arbitrary order along $A$ and
\[
\REa{i}{V/\pl}(*A)=\frac{\tilde{\Omega}^i(*A)}{g^*\Omega^1_{\pl}\wedge\tilde{\Omega}^{i-1}(*A)}
\]
By ~\ref{8abr01} we have
$\diss{j}{\REa{i}{V/\pl}(*A)}=0, j>0$. We have the following long
sequence
\begin{equation}
\A' :=\  0\rightarrow
\REa{0}{V/\pl}(*A)\stackrel{{d'}^0}{\rightarrow}
\REa{1}{V/\pl}(*A)\stackrel{{d'}^1}{\rightarrow}\cdots\stackrel{{d'}^{n-1}}{\rightarrow}
\REa{n}{V/\pl}(*A)\rightarrow 0
\end{equation}
Now we would like to verify the hypothesis of Lemma
~\ref{diagram} for $\A'$ and $\A$. The maps $h$ are inclusions.
The only non-trivial hypothesis is the isomorphism
~(\ref{18mai01}) in a point $p\in A$. Choose a Stein neighborhood
$U$ and a coordinate system $(z_1,z_2,\ldots,z_n,t)$ around
$p\in A$ such that in this system $p=0$, $A$ is given by $z_1=0$
and $L_{t_0}$ by $t=t_0$. We  have proved in ~\ref{erot} in Appendix
~\ref{gudiny} that in $U$
\[
\frac{Ker ({d'}^i)}{Im {({d'}^{i-1})}}=\frac{Ker ({d}^i)}{Im
{({d}^{i-1})}}=0, i\geq 2
\]
\[
\frac{Ker ({d'}^1)}{Im {({d'}^{0})}}=\frac{Ker ({d}^1)}{Im
{({d}^{0})}}=\{p(t)\frac{dz_1}{z_1}\mid p(t)\in \O_{\C,0}\}
\]
which proves the desired isomorphism in $p\in A$. The conclusion
is that:
\begin{equation}
\label{23ago}
 \tih{i}\cong \frac{Ker(d_*^i)}{Im(d_*^{i-1})}, 0\leq i\leq n-1,\
 {'\H}^n\cong \frac{g_*\REab{n}{V/\pl}}{Im(d_*^{n-1})}
\end{equation}
where
\begin{equation}
\label{28ago01}
 0\stackrel{d_*^{-1}}{\rightarrow}
g_*\REab{0}{V/\pl}\stackrel{d_*^0}{\rightarrow}
g_*\REab{1}{V/\pl}\stackrel{d_*^1}{\rightarrow}\cdots\stackrel{d_*^{n-1}}{\rightarrow}
g_*\REab{n}{V/\pl}
\rightarrow 0
\end{equation}

\begin{lem}
\label{17mai01}
Let $X$ be a paracompact Hausdorff space and
\begin{equation}
\label{iran}
 0\rightarrow F\stackrel{i}{\rightarrow}
F_0\stackrel{d^0}{\rightarrow} F_1\stackrel{d
^1}{\rightarrow}\cdots\stackrel{d^{n-1}}{\rightarrow} F_n
\end{equation}
an exact sequence of sheaves of Abelian groups. Let also  $Y$ be
another paracompact Hausdorff space and $g: X\rightarrow Y$ a
continuous map. Suppose that
\begin{equation}
\label{sher}
\diss{q}{F_p}=0,\ \forall\ q\geq 1, n>p\geq 0
\end{equation}
 Then
\[
\diss{p}{F}\cong\frac{Ker (d_*^{p})}{Im (d_*^{p-1})},\  0\leq p<n
\]
and there exists a natural inclusion $\diss{n}{F}\rightarrow
\frac{g_*F_n}{Im (d_*^{n-1})}$ such that we have
\[
\ses {\diss{n}{F}} {\frac{g_*F_n}{Im (d_*^{n-1})}}
{\frac{g_*F_n}{g_*Im (d^{n-1})}}
\]
where $d_*^p$'s define the sequence
\[
0\stackrel{d_*^{-1}}{\rightarrow} g_* F_0
\stackrel{d_*^0}{\rightarrow} g_*F_1\stackrel{d_*
^1}{\rightarrow}\cdots\stackrel{d_*^{n-1}}{\rightarrow} g_*F_n
\]
\end{lem}
\begin{proof}
 The proof  is a slight modification of Lemma 6.3 of
\cite{nar}. Put $Z_p=Ker(d^p),0\leq p<n $. The first statement is
trivial for $p=0$. Therefore let us prove the first statement for
$p\geq 1$. The exactness of ~(\ref{iran}) at $F_{p}$ gives us
\begin{equation}
\label{khosh} \ses{Z_{p-1}}{F_{p-1}}{Z_p},\ 1\leq p < n
\end{equation}
and we get the long exact sequence
\[
\cdots\rightarrow\diss{q}{F_{p-1}}\rightarrow
\diss{q}{Z_p}\rightarrow \diss{q+1}{Z_{p-1}}\rightarrow
\diss{q+1}{F_{p-1}}\rightarrow\cdots
\]
By ~(\ref{sher}) we conclude that
\[
\diss{q}{Z_p} \cong \diss{q+1}{Z_{p-1}}, 1\leq p<n, q\geq 1
\]
Since $F\cong Z_0$, we have
\begin{equation}
\label{tanha}
\diss{p}{F}\cong\diss{p-1}{Z_1}\cong\cdots\cong\diss{1}{Z_{p-1}},
1\leq p \leq n
\end{equation}
~(\ref{khosh}) gives us also
\[
g_*F_{p-1}\stackrel{d^{p-1}_*}{\rightarrow} g_*Z_p\rightarrow
\diss{1}{Z_{p-1}}\rightarrow 0,\ 1\leq p<n
\]
and thus
\[
\diss{1}{Z_{p-1}}\cong \frac{g_*Z_{p-1}}{Im (d^{p-1}_*)}=
\frac{Ker (d^{p}_*)}{Im (d^{p-1}_*)}, 0\leq p < n
\]
We have proved the first part of the lemma. Now let us prove the
second part. We have the short exact sequence
\[
\ses {g_* Im (d^{n-1})}{g_*F_n}{ \frac{g_*F_n}{g_*Im (d^{n-1})}}
\]
$Im (d_*^{n-1})$ is a subsheaf of both ${g_* Im (d^{n-1})}$ and
${g_* F_n}$ so we can rewrite the above exact sequence as:
\begin{equation}
\label{migooyam} \ses {\frac{g_* Im (d^{n-1})}{Im
d_*^{n-1}}}{\frac{g_*F_n}{Im (d_*^{n-1})} }{\frac{g_*F_n}{g_*Im
(d^{n-1})}}
\end{equation}
The short exact sequence
\[
\ses{Z_{n-1}}{F_{n-1}}{Im d^{n-1}}
 \]
 gives us
\[
g_*F_{n-1}\stackrel{d^{n-1}_*}{\rightarrow}
g_*Im(d^{n-1})\rightarrow \diss{1}{Z_{n-1}}\rightarrow 0
\]
Therefore by ~(\ref{tanha}) we have
\begin{equation}
\label{migooyam1} \diss{n}{F}=\frac{g_* Im (d^{n-1})}{Im
(d_*^{n-1})}
\end{equation}
Note that for this we do not need to have
$\diss{q}{Im(d^{n-1})}=0,\forall q\geq 1$. Now ~(\ref{migooyam})
and ~(\ref{migooyam1}) finish the proof. 
\end{proof}

Since $V-A\stackrel{\pi}{\cong} M-\RR$ and $f_0\circ\pi=g\mid_{V-A}$, we can
use the symbol  $f_0$ instead of $g\mid_{V-A}$. The following
sequence
\begin{equation}
 0\rightarrow f_{0}^*\O_{\pl}\stackrel{i}{\rightarrow}
\REab{0}{V/\pl}\stackrel{d^0}{\rightarrow}
\REab{1}{V/\pl}\stackrel{d^1}{\rightarrow}\cdots\stackrel{d^{n-1}}{\rightarrow}
\REab{n}{V/\pl} \hbox{ in } V-A
\end{equation}
is exact even in the critical points of $f_0$ (see \cite{bri}, Proposition 1.7, iii).
 We can apply Lemma ~\ref{17mai01} to the above
sequence and obtain
\[
R^if_{0*}{(f_0^*\O_\pl)}\cong\frac{Ker (d_*^i)}{Im (d_*^{i-1})}\cong
\tih{i} , \ i<n
\]
\[
\ses{R^nf_{0*}{(f_0^*\O_\pl)}} {\frac{g_* \REab{n}{V/\pl}}{Im
(d_*^{n-1})} } {\frac{g_* \REab{n}{V/\pl}}{g_*Im (d^{n-1})}}
\]
where $d_*^i$ is defined in ~(\ref{28ago01}).

Let ${'H}(p_i)$ be the Brieskorn module of a singularity
$p_i$ of $g$ (see the first paragraph of Appendix ~\ref{gudiny}). Define
\begin{equation}
\label{esmesag}
\CC_c=
    \left\{ \begin{array}{ll}
        0  & c \hbox{ is a regular value} \\
        \oplus_{i}{'H}(p_i) & \hbox{ $p_i$'s are the critical points
within $L_c$ }
        \\
       \end{array}
    \right.
\end{equation}
Each stalk $\CC_c$ is a free
$\O_{\pl,c}$-module of rank $\mu_c$. There is defined a natural
restriction map
\begin{equation}
\label{23ago01}
 \pi: {'\H^n}\rightarrow \CC
\end{equation}
We denote by $\CC'$ its image. Now fix a critical value $c\in C$. The 
stalk $\CC_c'$ is a
$\O_{\pl,c}$-submodule of $\CC_c$.
\begin{lem}
\label{birahe} We have
\[
\CC'=\frac{g_* \REab{n}{V/\pl}}{g_*Im (d^{n-1})}
\]
\end{lem}
\begin{proof}
By ~(\ref{23ago}) for $i=n$,
$\CC'$ is the image of $g_* \REa{n}{V/\pl}$ under the projection
$ g_* \REa{n}{V/\pl}\rightarrow \CC$. The kernel of this map is
exactly $g_*Im (d^{n-1})$ and so the proof is finished.
\end{proof}
We know that
\[
(R^if_{0*}{(f_0^*\O_\pl))\mid_{\pl-C}
}\cong(\diss{i}{\C_{V-A}}\otimes_{\C_\pl}\O_\pl)\mid_{\pl-C}=\tilde{\H}^i
\]
(see \cite{bri} p. 120). By lemma ~\ref{birahe} we conclude that
\begin{equation}
\label{akhar1} \tih{i}\cong\tilde{\H}^i \ i<n,\ 
{'\H}^n\mid_{\pl-C}\cong\tilde{\H}^n\mid_{\pl-C}
\end{equation}
\begin{equation}
\label{akhar2} \ses{\tilde{\H}^n}{{'\H}^n}{\CC'}
\end{equation}
~(\ref{akhar1}) and ~\ref{20ago01}  imply
that $\tih{i},\ i<n$ (resp. ${'\H}^n\mid_{\pl-C}$)
is a freely generated sheaf of rank $\beta_i$ (resp. $\beta_n$).
Now for $c\in C$, since the division of ${'\H}^n_c$ 
over the freely generated of rank $\beta_n-\mu_c$ $\O_{\pl,c}$-module
$\tilde{\H}^n_c$ (see ~\ref{20ago01}) is freely generated of rank 
$\mu_c$ (see ~\ref{kalisch}), we conclude
that ${'\H}^n_c$ is freely generated of rank $\beta_n$. 

Consider a continuous family $\{\delta_t\}_{t\in U}$  of
$i$-dimensional cycles in $g^{-1}(U)-A$ in such a way that
$\delta_t$ lies in $L_t$. For any $\omega\in
\Omega^i(U)$  the integral $\int_{\delta_t}\omega$ is
well-defined. Let $\gamma$ be a path in $U$ going around $t$
anti-clockwise and $\Gamma$ be the surface in $V$ formed by the
union of curves $\Gamma=\cup_{s\in\gamma}\delta_s$. 
With the above notation we have
\begin{equation}
\label{anha}
\int_{\delta_t}\omega=\frac{1}{2\pi i}\int_\Gamma
\frac{df\wedge\omega}{f-t}, \ \
\frac{d}{dt}\int_{\delta_t}\omega=\int_{\delta_t}\nabla_{\frac{\partial}
{\partial t}}\omega
\end{equation}
For the proof of above equalities see \cite{arn}.
By the second formula in ~(\ref{anha}) we can see that the flat sections
of $\nabla$ in $\H^i$ go to the flat sections of $\nabla$ in $\tilde\H^i$ by the isomorphism in ~(\ref{akhar1}) and we know that
this isomorphism is
obtained by restriction of $\omega\in\H^i(U)$ to the fibers $L_t, t\in U$.
 This implies that this isomorphism sends 
$(\H^i,\nabla)$ to $(\tilde\H^i,\nabla)$.
The proof
of Theorem ~\ref{khodaya} for $i=0,1,\ldots,n-1,'n$ is finished.

Now let us prove Theorem ~\ref{khodaya} for $i=n,''n$.
There is a natural inclusion $\tih{n}\subset {'\H}^n$. Let $U$ be
a small open disk in $\pl$ and $t$ a regular holomorphic function
in $U$. By Lemma ~\ref{dohafte} we have also the inclusion
\[
{'\H}^n\mid_U\stackrel{dt\wedge}{\rightarrow} {''\H}^n\mid_U
\]
We can see that
\[
\frac{{''\H}^n}{{'\H}^n}\mid_U\cong \REaa{n+1}{V/\pl}\mid_U
\]
and
$\frac{{'\H}^n}{\tih{n}}\stackrel{d(.)}{\rightarrow}\frac{{''\H}^n}{{'\H}^n}$ is
 an inclusion and so by ~\ref{kons} we conclude that 
$\H^n, {''\H}^n$
 are locally free sheaves of rank $\beta_n$. If $U\subset \pl-C$ then the above
inclusions are isomorphism of sheaves with connections.\qed

By the first part of Corollary ~\ref{18mar01} we know that
${\cal C}_c/{\cal C}_c,\ c\in C$ is a vector space of dimension less 
than $\mu_c$. I believe that it is zero. 
\appendix
\section{Complex Geometry}
\label{17mar01}
\def\T{{\cal T}}
\def\L{{\cal L}}
\newcommand{\dis}[2]{R^{#1}f_{*}#2}
\newcommand{\disg}[2]{R^{#1}g_*#2}

In this appendix we will give all preliminaries
in complex analysis and complex geometry used throughout the
article. I did not  find a book in the literature of complex
analysis containing all of these preliminaries and so I have
collected them in this appendix.

In what follows by an analytic sheaf over an analytic variety $V$
we mean a $\O_V$-module sheaf. For a given analytic sheaf $\s$
over $V$, when we write $x\in \s$ we mean that $x$ is a
holomorphic section of $\s$ over some open neighborhood in $V$ or
it is an element of some stalk of $\s$; being clear from the text
which we mean.

{\bf Direct Limit Sheaves:}
Let $\{\s_i\}_i$ be a direct system of sheaves i.e.,
\[
\s_0\rightarrow\s_1\rightarrow\cdots\s_i\rightarrow\cdots
\]
If there is no confusion we write simply $\{\s_i\}$.
We define the direct limit of the system, say $\lim_{i\to\infty} \s_i$, to be
the sheaf associated to the presheaf $U\rightarrow \lim_{i\to\infty} \s_i(U)$.
There are defined natural maps $\s_i\rightarrow \lim_{i\to\infty}\s_i$.

Let $\s$ be another analytic sheaf and $\{\s_i\rightarrow \s\}$ a
collection of compatible analytic homomorphisms. Then there is a
unique map $\lim_{i\to\infty} \s_i\rightarrow \s$ such that for each $i$, the
original map $\s_i\rightarrow \s$ is obtained by composing the
maps $\s_i\rightarrow \lim_{i\to\infty} \s_i\rightarrow \s$.
\begin{app}
\label{del} Let $\{\s_i\}$ be a direct system of sheaves and
$\{\s_i\rightarrow \s\}$ a collection of compatible maps. Then
$\lim_{i\to\infty} \s_i\rightarrow \s$ is an isomorphism if and only if
\begin{enumerate}
\item
For any $x\in\s$ there exist $i\in \N$ and $x_i\in \s_i$ such that
$x_i\rightarrow x$;
\item
If there exist $i_0\in \N$ and  a sequence $x_{i_0}\rightarrow x_{i_0+1}
\rightarrow\cdots,\ x_i\in \s_i$ such that $x_i\rightarrow 0\in \s$
then there exists $i_1\geq i_0$ such that for all $i\geq i_1$ we have $x_i=0$.
\end{enumerate}
\end{app}
\begin{proof}
The first statement implies the surjectivity and the second one implies the
 injectivity of $\lim_{i\to\infty}\s_i\rightarrow \s$.
\end{proof}

 Using the above proposition we can check the following simple
 facts:
 \begin{app}\rm
 \label{A1}
The short exact sequence
\[
\ses{\{\L_i\}}{\{\s_i\}}{\{\T_i\}}
\]
gives
\[
\ses{\lim_{i\to\infty}\L_i}{\lim_{i\to\infty}\s_i}{\lim_{i\to\infty}\T_i}
\]
\end{app}
\begin{app}\rm
\label{A2}
 For a collection of compatible maps
$\{\s_i\}\rightarrow\{\T_i\}$ we have
\[
\lim_{i\to\infty} Ker(\{\s_i\}\rightarrow \{\T_i\})=Ker(\lim_{i\to\infty} \s_i \rightarrow
\lim_{i\to\infty} \T_i)
\]
\[
\lim_{i\to\infty} Im(\{\s_i\}\rightarrow \{\T_i\})=Im(\lim_{i\to\infty} \s_i \rightarrow \lim_{i\to\infty} \T_i)
\]
\end{app}
One of the important properties of the direct limit sheaf is:
\begin{app}\rm
\label{A3}
Let $\{\s_i\}$ be a direct system of sheaves over $V$.
If $V$ is compact then
\[
H^\mu(V, \lim_{i\to\infty}\s_i)=\lim_{i\to\infty} H^\mu(V, \s_i), 
\mu=0,1,2,\ldots
\]
\end{app}
\begin{proof}
The trick of the proof is that for a finite covering ${\cal U}$ of $V$
with Stein open sets every $\alpha\in H^\mu({\cal U}, \lim_{i\to\infty} \s_i)$
($Z^\mu({\cal U}, \lim_{i\to\infty} \s_i)$ or $B^\mu({\cal U}, \lim_{i\to\infty} \s_i)$) is represented
by a finite number of sections. This enables us to check the properties 1
 and 2 of Proposition ~\ref{del}.
\end{proof}

{\bf Sheaves with Pole Divisors:} Let $\s$ be an analytic sheaf
over an analytic compact variety $V$ and $D$ a divisor in $V$
which does not intersect the singular locus of $V$. By $\s(kD)$
we denote the sheaf of meromorphic sections of $\s$ with poles of
multiplicity at most $k$ along $D$. Also, $\s(*D)=\lim_{k\to\infty} \s(kD)$
denotes the sheaf of meromorphic sections of $\s$ with poles
of arbitrary order along $D$. We list some natural properties of sheaves with poles.

\begin{app}\rm
\label{A4}
The short exact sequence of analytic sheaves $0\rightarrow \L\rightarrow \s\rightarrow \T\rightarrow 0$ gives us the short exact sequence $0\rightarrow \L(*D)\rightarrow \s(*D)\rightarrow \T(*D)\rightarrow 0$. In particular if $\L$ is a subsheaf of $\s$ then $\frac{\s}{\T}(*D)=\frac{\s(*D)}{\L(*D)}$.
\end{app}
\begin{app}\rm
\label{A5}
The analytic homomorphism of sheaves $d:\s\rightarrow \T$ induces a natural analytic homomorphism $d_D:\s(*D)\rightarrow \T(*D)$ and
\[
ker(d)(*D)=ker(d_D),\ Im(d)(*D)=Im(d_D)
\]
\end{app}
\begin{app}\rm
\label{A6}
Let $D$ be a divisor in $V$. We have
\[
(\lim_{i\to\infty} \s_i)(*D)=\lim_{i\to\infty} \s_i(*D)
\]
\end{app}
\begin{app}\rm
\label{A7}
$(\s\otimes_{\O_V}\T)(*D)=\s\otimes_{\O_V}\T(*D)=\s(*D)\otimes_{\O_V}\T$
\end{app}
\begin{app}\rm
\label{A8} If $\s$ is coherent then $\s(kD)$ is also coherent.
Moreover if $V$ is a compact manifold and $D$ is a positive
divisor then there exists an integer $k_0$ such that
\[
H^\mu(V,\s(kD))=0, \ k\geq k_0,\ \mu \geq 1
\]
Using A.4 and $\s(*D)=\lim_{k\rightarrow \infty}\s(kD)$ we have
\[
H^\mu(V,\s(*D))=0, \ \mu \geq 1
\]
\end{app}
\begin{app}\rm
\label{A9}
(Vanishing theorem for limit sheaves)
Let $\{\s_i\}$ be a direct system of coherent sheaves and $D$ a positive divisor in $V$.
 If $V$ is a compact manifold then
\[
H^\mu(V,\lim_{i\to\infty}\s_i(*D))=0,\ \mu\geq 1
\]
\end{app}
\begin{proof}
This is a direct consequence of ~\ref{A3} and ~\ref{A8}.
\end{proof}
Let $d:S\rightarrow\T$ be an analytic map between 
two coherent sheaves on a complex manifold $V$, $D$ a positive divisor
in $V$ and $H^0(d_D):H^0(V, S(*D))\rightarrow H^0(V, \T(*D))$. 
\begin{app}
\label{shayad}
We have $H^0(V, ker(d_D))=ker(H^0(d_D)),\ H^0(V, Im(d_D))=Im(H^0(d_D))$. \end{app}
The coherence of the sheaves and the positivity of the divisor is strongly
used in the second equality.

{\bf Direct Image Sheaves:} The first lines of this paragraph can
be found in  Chapter 1 Section 4.7 of \cite{gra2}. 
Let $f:X\rightarrow Y$ be a holomorphic map between the analytic
varieties $X$ and $Y$ and $\s$ an analytic sheaf on $X$. For any
open Stein subset $U$ of $Y$ we can associate the $\O(U)$-module
$H^i(f^{-1}(U), \s)$. There are canonical restriction maps and we
have an analytic presheaf on $Y$ defined on all open Stein
subsets of $Y$. The associated analytic sheaf on $Y$ is called
the $i$-th direct image of $\s$ and is denoted by $\dis{i}{\s}$.
Every short exact sequence  $\ses{\L}{\s}{\T}$ of analytic
sheaves over $X$ induces a long exact cohomology sequences
\[
0\rightarrow \dis{0}{\L}\rightarrow \dis{0}{\s}
\rightarrow\dis{0}{\T}\rightarrow\dis{1}{\L}\rightarrow\dis{1}{\s}
\rightarrow\cdots
\]
The following fact says that the functor $\dis{i}{}$ and $\lim$ commute:
\begin{app}\rm
Let $f:X\rightarrow Y$ be a holomorphic map between the analytic varieties
$X$ and $Y$ and $\{\s_k\}$ a direct system of  analytic sheaves over $X$. Then
\[
\lim_{k\to\infty} \dis{i}{\s_k}=\dis{i}\lim_{k\to\infty} \s_k
\]
\end{app}
 The Grauert direct image theorem says when the direct image
sheaf $\dis{i}{\s}$ is coherent:
\begin{app}\rm
(Grauert direct image theorem)
\label{gdit}
Let $f:X\rightarrow Y$ be a proper holomorphic map between the analytic varieties
$X$ and $Y$ and $\s$ a coherent analytic sheaf over $X$. Then for any $i\geq 0$
the $i$-th direct image $\dis{i}{\s}$ is a coherent analytic sheaf over $Y$.
\end{app}
Let $X'$ be an analytic subvariety of $X$ and $\s$ an analytic
sheaf over $X$. By structural restriction of $\s$ to $X'$ we mean
$\s\mid_{X'}=\frac{\s}{\M.\s}$, where $\M$ is the sheaf of
holomorphic functions vanishing on $X'$. If $\s$ is a coherent
$\O_X$-module sheaf then $\s\mid_{X'}$ is a coherent
$\O_{X'}$-module sheaf. This restriction is different with the
 sheaf theorical restriction. In what follows all restrictions
 we consider are structural except in mentioned cases.

Let $g:V\rightarrow \pl$ be the holomorphic function introduced in
in the first section, $c$ a point in $\pl$ and $\s$ an analytic
sheaf on $V$.
\begin{app}
\label{A10.5}
 Let $\s$ be an analytic sheaf on $V$. Then
\[
(\diss{i}{\s})(*D)=\diss{i}(\s(*g^{-1}(D)))
\]
where $D=\{p\}$.
\end{app}
The above proposition in general may not be true 
(for instance when $g$ has multiplicity along $g^{-1}(D)$).
\\

 We define $\s_c=\s\mid_{V_c}$ to be
 the restriction of $\s$ to the fiber $V_c=g^{-1}(c)$.
 The following natural function is well-defined:
\[
g_{c,i}: \disg{i}{\s}\mid_c \rightarrow H^i(V_c,\s_c)
\]
\begin{app}\rm
\label{mahetaban} The map $g_{c,i},\ i\geq 0$ is injective.
\end{app}
 \begin{proof}
Suppose that for an $\alpha\in \disg{i}{\s}\mid_c$ we have
$g_{c,i}(\alpha)=0$. For a Stein covering $\U$ of $V_U$, $\alpha$
is represented by an element $\alpha\in H^i(\U,\s)$, where $U$ is
a small open disk around $c$. $g_{c,i}(\alpha)=0$ means that the
restriction of $\alpha$ to $V_c$ is zero. In other words there
exists a $\beta\in C^{i-1}(\U\cap V_c, \s\mid_{V_c})$ such that
$\alpha=\partial \beta$. Since $\U$ is a Stein covering, taking
$U$ smaller if it is necessary we can represent  $\beta$ as an
element of $C^{i-1}(\U, \s)$ (by extending $\beta$). Now
$\alpha-\partial\beta\mid_{V_c}=0$ and so
$\alpha-\partial\beta=h\circ g . \gamma$ for some $\gamma\in
C^{i}(\U, \s)$, where $h$ is a holomorphic regular function on $U$
vanishing on $c$ (here we have used this fact that the
multiplicity of $g$ along each irreducible component of $V_c$ is
one). Therefore $\alpha$ is zero in $\disg{i}{\s}\mid_c$. 
\end{proof}

The map $g_{c,i}$ need not to be surjective. The obstruction  to
the surjectivity of $g_{c,i}$ is an element
$\alpha\in\disg{i+1}{\s}_c$ with $supp(\alpha)=\{c\}$. Therefore
if $\disg{i+1}{\s}$ is freely generated then $g_{c,i}$ is an
isomorphism (For more information see \cite{gra2} p. 209).

Our main Theorem in this paragraph which is used frequently in
the article is the following:
\begin{app}\rm
\label{8abr01} (Variational vanishing theorem) Let $g:V\rightarrow
\pl$ be as before and $\s$ a coherent sheaf on $V$. Let also $A$
be the blow-up divisor in $V$. Then
\[
\disg{i}{\s(*A)}=0, \ i\geq 1
\]
\end{app}
\begin{proof}
The main property of $A$ is that it  its intersection with
each fiber $V_c$ is positive in $V_c$. Fix a regular value
$c\in\pl$. Since $A_c=A\cap V_c$ is positive in $V_c$, there
exists a natural number $k_0$ such that
\[
H^i(V_c, \s_c(kA_c))=0,\ k\geq k_0
\]
This and ~\ref{mahetaban}  imply that
$\disg{i}{\s(kA)}\mid_c=0$. By Grauert direct image theorem
$\disg{i}{\s(kA)}$ is coherent, therefore $\disg{i}{\s(kA)}$ is
the zero sheaf in a neighborhood of $c$. Now in this neighborhood
we have
\[
\disg{i}{\s(*A)}= \disg{i}{\lim_{k\to\infty} \s(kA)}=\lim_{k\to\infty}
\disg{i}{\s(kA)}=\lim_{k\to\infty} 0=0
\]
Until now we have proved that $supp(\disg{i}{\s(*A)})\subset C$. If
we had some type of Kodaira vanishing theorem for a singular
variety $V_c,c\in C$ then the proof was complete. But I do not
know such a theorem and so I use the following trick: Let $b$ be a
regular value in $\pl$. Since $\disg{i}{\s(*A)}$ is a discrete
sheaf, we have
\[
H^0(\pl,\disg{i}{\s(*A)}(*b))=\cup_{c\in C}\disg{i}{\s(*A)}_c
\]
By ~\ref{A10.5} we have $\disg{i}{\s(*A)}(b)=\disg{i}{\s(*A+*V_b)}$ 
and so
$H^0(\pl,\dis{i}{\s(*A)}(*b))=H^0(\pl,\disg{i}{\s(*A+*V_b)}
)\subset H^i(V,\s(*A+ *V_b))$. $A \cup V_b$ is the pullback of
$M_b$ by the blow up map $\pi:V\rightarrow M$ and $M_b$ is a
hyperplane section of $M$. Therefore $A \cup V_b$ is a positive
divisor and
\[
H^i(V,\s(*A +*V_b))=0
\]
We conclude that $\disg{i}{\s(*A)}_c=0$ for $c\in C$ which is the
desired. 
\end{proof}
{\bf Atiyah-Hodge type Theorems:} Let $V$ be a projective
manifold of dimension $n$ and $A$ a submanifold of $V$ of
codimension one. Denote by $\Omega^i(*A)$ the sheaf of meromorphic
$i$-forms in $V$ with poles of arbitrary order along $A$.
 We have the following, not necessarily exact, sequence
\[
0\rightarrow \C\stackrel{}{\rightarrow}H^0(V,\Omega^0(*A))
\stackrel{d^0}{\rightarrow} H^0(V, \Omega^1(*A))\stackrel{d^1}{\rightarrow}\cdots
\stackrel{d^{i-1}}{\rightarrow}H^0(V, \Omega^i(*A))\stackrel{d^{i}}{\rightarrow}\cdots
\]
We form the cohomology groups

\[
\tilde{H}^i=\frac{Ker(d^i)}{Im(d^{i-1})},\ i\geq 0
\]
\begin{app}(Atiyah-Hodge Theorem \cite{nar})
\label{hodge}
Suppose that $A$ is positive in $V$. Then
there are natural isomorphisms
\[
H^i(V-A,\C)\cong \tilde{H}^i
\]
\end{app}
Roughly speaking, this theorem says
that every cohomology class in $H^i(V-A,\C)$ is represented by a
closed meromorphic i-form in $V$ with poles along $V$.
\section{Some topological facts}
\label{alman}
All homologies considered in this appendix are with rational
coefficients. Recall the notations ~(\ref{22}),~(\ref{11}). 
Let $c\in\pl$, $U$ a small open disk with center $c$
and $b$ a regular point in the boundary of $U$. We denote by $D$
the closure of $U$ in $\pl$. Let also $\{p_i\mid i=1,2,\ldots,
k\}$ be the singularities within $L_c$. To each $p_i$ we can
associate a set of distinguished vanishing cycles
$\{\delta_{ij}\mid j=1,2,\ldots, l_k\}$ in $H_{n}(L_b)$ (see
\cite{arn}). Let also $\mu_c$ denote the sum of Milnor numbers of 
singularities 
within $L_c$. If $c$ is a regular value of $f_0$ then $\mu_c=0$.
\begin{app}
\label{20ago01}
 We have 1. $L_c$ is a deformation retract of $L_D$ 2. $H_{n+1}(L_D)=0$ 
3. $H_i(L_D,L_b)=0$ for $0\leq i \leq n$ and $H_{n+1}(L_D,L_b)$ is
freely generated of rank $\mu_c$
4. $L_U$ is a Stein manifold
5. There is no linear relation between $\delta_{ij}$'s
6. $\tilde{\H}^i_c, 0\leq i\leq n-1$ is a freely generated
$\O_{{\pl},c}$-module of rank $\beta_i$ and $\tilde{\H}^n_c$ is a
freely generated $\O_{{\pl},c}$-module of rank $\beta_n-\mu_c$.
\end{app} 
\begin{proof}
 Let us prove the first part. Since out of $c$ the map $g$
 is a $C^\infty$ fiber bundle, by homotopy covering theorem
 (see 14, 11.3,\cite{ste}) we can take $U$ smaller if it is
 necessary.
 Let $B_i,i=1,2,\ldots,k$ be an open ball with center $p_i$
 whose boundary is transverse to $L_t, t\in D$.
 $f:(L_D-\cup_i B_i, \partial (L_D-\cup_i B_i) )\rightarrow D$ is a $C^\infty$ fibration.
 Therefore $L_D$ can be retracted to
 $L_c\cup\cup_i (L_D\cap B_i)$. Now by an argument stated in
 \cite{arn} p.32 we know that $L_c\cap B_i$ is a deformation
 retract of $L_D\cap B_i$ and so $L_c$ is a deformation retract of
 $L_D$.
 
Let us prove the second part. Let $\delta$ be an $(n+1)$-cycle in
$L_D$. Taking another cycle
 in the homological class of $\delta$ we can assume that $\delta$
 does not pass through $p_i$'s. This time we take the ball
 $B_i$ in such a way that it does not intersect $\delta$. Let
 $D'$ be another small closed disk inside $D$ with center $c$ such that
 $L_t,t\in D'$ is transverse to $\partial B_i$. $L_{D'}$ is a deformation
 retract of $L_{D}$ and $L_b\cup\cup_i (L_{D'}\cap B_i)$  is a deformation retract of
 $L_{D'}$, where $b$ is a regular
 value in the boundary of $D'$. Therefore $\delta$ is
 homologous to an $(n+1)$-cycle in $L_b$. But $L_b$ is a Stein
 manifold of dimension $n$ and so
 $H_{n+1}(L_b)=0$. We conclude that $H_{n+1}(L_U)=0$.
 
 The proof of the third part is the same as (5.4.1) of
 \cite{lam}. Instead of (5.5.9)\cite{lam} we use  a
similar statement for an arbitrary isolated singularity 
(see \cite{arn}).
 
$L_U$ has no non discrete compact analytic set because
 $L_U=\cup_{t\in U}L_t$ and each $L_t$ is a Stein analytic
space. Now let us prove that $L_U$ is holomorphically convex. To
see this fact let $p=\infty\not\in U$ and $t\in U$. We can
consider $L_U$ as a subset of the Stein manifold
$M-M_{p}=L_{\C}$.
 Every holomorphic function in $L_t$ extends to $M-M_{p}$ and hence to
$L_U$. Knowing this and the fact that each $L_t,t\in U$ is Stein,
we can easily check that $L_U$ is holomorphically convex.

Now let us prove the fifth part.  Writing the long exact sequence
of the pair $(L_D,L_b)$ we have:
\begin{equation}
\label{mordam}
 \cdots\rightarrow H_{n+1}(L_D)\rightarrow
H_{n+1}(L_D,L_b)\rightarrow H_n(L_b)\rightarrow
H_n(L_D)\rightarrow 0
\end{equation}
$H_{n+1}(L_D)=0$ and the vanishing cycles $\delta_{ij}$ are images
of a basis of $H_{n+1}(L_D,L_b)$ under the boundary map. Therefore
there does not exist any linear relation between $\delta_{ij}$'s.

Now let us prove the last part. Since $L_c$ is a deformation
retract of $L_D$, $\tilde{\H}^i_c, i\leq n$ is freely generated
of rank $dimH_i(L_c)$. By 1,2,3 and the long exact sequence of
the pair $(L_D,L_b)$, we have $dim(H_i(L_c))=\beta_i,i\leq n-1$
and $dim(H_n(L_c))= \beta_n-\mu_c$.
\end{proof}
Both the inclusions $L_b\subset L_D ,b\in D$ and $L_c\subset L_D$
induce isomorphisms in $i$-th homologies, where $i\leq n-1$ and
$i=n$ if $c$ is a regular value. This means that we have a natural
$i$-th homology bundle, and hence $i$-th cohomology bundle over
$\pl$, for $i\leq n-1$ (for $i=n$ over $\pl-C$).
\section{Local Brieskorn modules}
\label{gudiny}
Let $f:(\C^{n+1},0)\rightarrow(\C,0)$  be a germ of a holomorphic function with
an isolated critical point at $0$. The
Brieskorn module
\[
'H(0)=\frac{\Omega^{n}}{df\wedge\Omega^{n-1}+d\Omega^{n-1}}
\]
is a freely generated $\O_{\C,0}$-module of rank $\mu$ (see
\cite{bri} and \cite{seb}), where $\Omega^i$ is the set of
$i$-forms in $(\C^{n+1},0)$ and $\mu$ is the Milnor number of $f$.
Let $t$ be a coordinate system in $(\C,0)$.
\begin{app}
\label{puf}
 Suppose that the restriction of $\{\omega_i\in \Omega^n\mid 1\leq i\leq \mu\}$
to a fiber $f^{-1}(t),t\in(\C,0)-\{0\}$ generates its cohomology group. 
Then for all
 $\omega\in \Omega^n$ there exists a natural number $h$ such that
 $t^h\omega$ belongs to the $\O_{\C,0}$-module generated by $\omega_i$'s in
 ${'H}(0)$.
\end{app}
\begin{proof}
Let $\{\delta_j(t)\mid 1\leq j\leq \mu\}$ be a basis of
vanishing cycles in $H_n(f^{-1}(t), \Z)$. Define the matrices
$A=[\int_{\delta_j(t)}\omega_i]_{\mu\times\mu}$ and
$B=[\int_{\delta_j(t)}\omega]_{\mu\times 1}$. Define
\begin{equation}
\label{ahhhhh}
 P:=A^{-1}B=\frac{adj(A).B}{det(A)}
\end{equation}
where $adj(A)$ is the adjoint of $A$. If we change the basis of
$H_n(f^{-1}(t),\Z)$ and $C$ is the matrix of this change then $A$
changes to $C.A$ and $B$ to $C.B$, therefore
$P=(C.A)^{-1}C.B=A^{-1}B$ does not change, particularly when $C$
is the monodromy operator obtained by turning around $0$. We
conclude that $P$ is a one valued holomorphic function in
$(\C,0)-\{0\}$. ~(\ref{ahhhhh}) implies that$P=\frac{P'}{t^h}$,
where $h$ is a natural number and $P'=[p_j],p_j\in\O_{\C,0}$. Now
\[
\int_{\delta_{j}(t)}(t^h\omega-\sum_i p_i\omega_i)=0 \ \forall
\delta_j(t)
\]
A basis of the freely generated $\O_{\C,0}$-module $'H(0)$
generates $H^n(f^{-1}(t),\C),\forall t\in(\C,0)-\{0\}$. Therefore
$t^h\omega-\sum_i p_i\omega_i$ is zero in ${'H}(0)$.
 \end{proof}
Recall 
~(\ref{esmesag}) and ~(\ref{23ago01}). The main proposition in this appendix 
is:
\begin{app}
\label{kalisch}
$\CC'_c$ is a free $\O_{\pl,c}$-module of rank $\mu_c$.
\end{app}
Its proof consists of various steps. 
\begin{app}
\label{9ago01} 
$\frac{\CC_c}{\CC'_c}$ is a finite dimensional
vector space.
\end{app}
\begin{proof}
Let $a$ be a regular point in $U$. Since  $A\cap V_a$ is a hyperplane
section of $V_a$ and $\delta_{ij},\  
i=1,2,\ldots,k\
j=1,2,\ldots,l_k$ are linearly independent,
 by Atiyah-Hodge theorem there are
meromorphic $n$-forms $\omega_{ij}$ in $V_a$ with poles along $A\cap V_a$ such that
$det[\int_{\delta_{ij}}\omega_{ij}]_{\mu_c\times\mu_c}\not= 0$

Consider the sheaf $\s$ of holomorphic $n$-forms in $V_U$ which
are zero restricted to $V_a$. $\s$ is a coherent sheaf and so by
 ~\ref{8abr01} $H^1(V_U,\s(*A))=0$, where $\s(*A)$ is the
sheaf of meromorphic sections of $\s$ with poles of arbitrary order
along $A$. This
implies that each $\omega_{ij}$ extends to $V_U$ as a meromorphic
$n$-form  with poles along $A$. We use the same notations for the extended ones.
We conclude that the restriction
of the $n$-forms $\omega_{ij}$ to a regular fiber of a singularity 
$g:(V,p_i)\rightarrow (\pl,c)$ generate its $n$-th cohomology
group. This and ~\ref{puf}  imply that for every
$\omega\in \CC_c/\CC'_c$ there exists a natural number $h$ such
that $(t-c)^h\omega=0$. Let $\Omega=\{\omega_i\mid i=1,2,\ldots,
\mu_c\}$ freely generate $\CC_c$ and $h$ be the minimum number
such that $(t-c)^h\Omega=0$ in $\CC_c/\CC'_c$. Now $\cup_{0\leq
i\leq h-1}(t-c)^i\Omega$ generates $\CC_c/\CC'_c$ as a vector
space. 
\end{proof}
\begin{app}
\label{kons}
If $S$ is a free ${\cal O}_{\C,0}$ module of finite rank $k$ and if 
$R$ is a submodule, then $R$ is free of rank $l\leq k$; one has $l=k$ if and 
only if $\dim_\C S/K<\infty$.  
\end{app}
$\O_{\C,0}$ is a principal ideal domain and 
the proof follows form the structure theroy of modules of principal ideal 
domains.
{\it Proof of ~\ref{kalisch}.}
 We know that $\CC_c$ is a free $\O_{\pl,c}$-module of rank
$\mu_c$. Also by ~\ref{9ago01}, $\CC_c/\CC'_c$ is a
finite dimensional vector space. Therefore we can apply
~\ref{kons} and conclude the theorem.\qed

I believe that $\CC'_c=\CC_c$. But the methods in this article are not
sufficient to prove this stronger result.

Let $A=\{(z_1,z_2,\ldots,z_n)\in(\C^n,0)\mid z_1=0\}$ and
$\Omega^i$  be the set of holomorphic $i$-forms in $(\C^n, 0)-A$.
We have the complex
\[
 0\rightarrow \Omega^0 \stackrel{{d}^0}{\rightarrow}
\Omega^1\stackrel{{d}^1}{\rightarrow}\cdots\stackrel{{d}^{n-1}}{\rightarrow}
\Omega^n\stackrel{d^n}{\rightarrow} 0
\]
and so we define $H^i=\frac{Ker(d^i)}{Im(d^{i-1})}$. Let also
$\Omega^i_*$ be the subset of $\Omega^i$ containing the $i$-forms
with poles of arbitrary order along $A$. In the same way we can define
$H^i_*$.
\begin{app}
\label{erot} we have $H^i=H^i_*=0, i\geq 2$ and $H^1=H^1_*=L$,
where $ L=\{p\frac{dz_1}{z_1}\mid p\in \O_{\C^k,0}\}$
\end{app}
\begin{proof}
 The proof is completely formal, for instance see \cite{gun}
Theorem 3E. We only prove the proposition for $H^i$. For the
other the argument is similar.

Fix the $i$-form $\omega$ with $d\omega=0$. We want to prove that
$\omega=d\eta$ (up to $L$ if $i=1$), where $\eta\in\Omega^{i-1}$.
Let $k$ be the least integer such that the representation of
$\omega$ contains only $dz_1,dz_2,\ldots,dz_k$ (we have $i\leq
k$). The proof is by induction on $k$. The case $k=0$ is trivial.
We write $\omega=dz_k\wedge \alpha+\beta$, where $\alpha$ and
$\beta$ are differential forms that involve only
$dz_1,dz_2,\ldots,dz_{k-1}$. Since $ -dz_k\wedge d\alpha
+d\beta=0$, the coefficients of $\alpha$ and $\beta$ do not
depend on $z_{k+1},\ldots,z_n$. If $k>1$ then we can write any
coefficient of $\alpha$, say $f$, as $f=\frac{dg}{dz_k}$ and if
$k=1$ as $f=\frac{p}{z_1}+\frac{dg}{dz_1}$, where $g\in \Omega^0$
and $p\in\O_{\C^k,0}$.
 Let
$\gamma$ be the differential form obtained from $\alpha$ by
replacing each coefficient $f$ by the corresponding coefficient
$g$. Then $d\gamma=\delta+ dz_k\wedge\alpha$ (if $k=1$ then up to
$L$), where $\delta$ is a differential form involving only
$dz_1,\ldots,dz_{k-1}$. Next set $\theta=\delta-\beta$. We have
$d\theta=0$ and so by induction $\theta=d\eta$ (if $i=1$ then up
to $L$). Now $\omega=d(\gamma-\eta)$ (if $i=1$ then up to $L$).
\end{proof}

\end{document}